\documentclass[journal]{IEEEtran}

\pdfoutput=1
\usepackage{enumitem}
\usepackage{amsmath, amssymb}
\usepackage[pdftex]{graphicx}
\graphicspath{{./Figures/}} 
\usepackage{epstopdf}
\usepackage{float}
\usepackage{subcaption} 
\usepackage{booktabs}
\usepackage{threeparttable}
\usepackage{tablefootnote}
\usepackage{multirow}
\usepackage{dsfont}
\usepackage{amsthm} 

    \newtheorem{remark}{Remark}

    \newtheorem{proposition}{Proposition}

\usepackage{enumitem}
\usepackage{threeparttable}
\usepackage{tablefootnote}
\usepackage[hang,flushmargin]{footmisc}

\usepackage{url}

\usepackage[ruled, lined, ,longend, linesnumbered]{algorithm2e}
\usepackage{algpseudocode}
\SetKwInOut{Input}{Input}
\SetKw{KwBy}{by}
\SetKw{Kwor}{or}
\SetKw{Kwand}{and}
\SetKw{KwEndFor}{end for}

\makeatletter
\newcommand\notsotiny{\@setfontsize\notsotiny\@vipt\@viipt}
\makeatother
\usepackage{fancyhdr}

\def\set#1#2{\{ \, #1 :#2\,\}} 
\newcommand{\cc}[1]{{\mathcal{#1}}} 
\def\R{{\mathbb{R}}} 
\def\Z{{\mathbb{I}}} 
\def\inR#1{\in\R^{#1}} 
\def\vv#1{{ \rm \bf{#1}}} 
\def\Dp#1{\mathbb{D}_{+}^{#1}} 
\def\Sp#1{\mathbb{S}_{+}^{#1}} 
\def\Ssp#1{\mathbb{S}^{#1}} 
\newcommand{\T}{^\top} 
\def\fracg#1#2{{\displaystyle{\frac{#1}{#2}}}} 
\def\Sum#1#2{\sum\limits_{#1}^{#2}} 
\DeclareMathSymbol{\shortminus}{\mathbin}{AMSa}{"39} 

\def\bmat#1{\left[\begin{array}{#1}} 
\def\emat{\end{array}\right]} 
\newcommand {\bsis} {\left\{ \begin{array} }
\newcommand {\esis} {\end{array}\right.}
\newcommand{\bRlist}{\renewcommand{\labelenumi}{(\roman{enumi})} \begin{enumerate}} 
\newcommand{\eRlist}{\end{enumerate} \renewcommand{\labelenumi}{\arabic{enumi}}} 

\def\E{\mathcal{E}} 

\def\ta{\tilde a}
\def\tx{\tilde x}
\def\tu{\tilde u}
\def\zi{z_\circ}
\def\vi{v_\circ}
\def\lbvi{\underline{v}_\circ}
\def\ubvi{\overline{v}_\circ}
\def\lambdai{\lambda_\circ}
\def\zN{z_f}
\def\vN{v_f}
\def\lambdaN{\lambda_f}

\begin{document}
\pagestyle{fancy}

\title{\LARGE A sparse ADMM-based solver for linear MPC subject to terminal quadratic constraint%
}

\author{
    Pablo~Krupa$^{\star}$,~Rim~Jaouani$^{\dagger}$,~Daniel~Limon$^\dagger$,~Teodoro~Alamo$^\dagger$%
    \thanks{$^\star$ Gran Sasso Science Institute (GSSI), L'Aquila, Italy. Corresponding author. E-mail: \texttt{pablo.krupa@gssi.it}.}%
    \thanks{$^\dagger$ Systems Engineering and Automation, Universidad de Sevilla, Sevilla, Spain. E-mails: \texttt{rjaouani@us.es}, \texttt{dml@us.es}, \texttt{talamo@us.es}.}%
    \thanks{This work has been partially funded by grant PID2022-141159OB-I00 funded by MCIN/AEI/10.13039/501100011033 and by ERDF/EU, and by grant PDC2021-121120-C21 funded by MCIN/AEI/10.13039/501100011033 and by the ``European Union NextGenerationEU/PRTR".}
    \thanks{Pablo Krupa acknowledges the support of the MUR-PRO3 project on Software Quality and the MUR-PRIN project DREAM (20228FT78M).}
}

\maketitle
\thispagestyle{fancy}

\begin{abstract}Model Predictive Control (MPC) typically includes a terminal constraint to guarantee stability of the closed-loop system under nominal conditions.
In linear MPC this constraint is generally taken on a polyhedral set, leading to a quadratic optimization problem.
However, the use of an ellipsoidal terminal constraint may be desirable, leading to an optimization problem with a quadratic constraint.
In this case, the optimization problem can be solved using Second Order Cone (SOC) programming solvers, since the quadratic constraint can be posed as a SOC constraint, at the expense of adding additional slack variables and possibly compromising the simple structure of the solver ingredients.
In this paper we present a sparse solver for linear MPC subject to a terminal ellipsoidal constraint based on the alternating direction method of multipliers algorithm in which we directly deal with the quadratic constraints without having to resort to the use of a SOC constraint nor the inclusion of additional decision variables.
The solver is suitable for its use in embedded systems, since it is sparse, has a small memory footprint and requires no external libraries.
We compare its performance against other approaches from the literature.
\end{abstract}


\subsubsection*{Notation}
We denote by $\Ssp{n}$, $\Sp{n}$ and $\Dp{n}$ the spaces of positive semi-definite, positive definite, and diagonal matrices with positive entries of dimension $n \times n$, respectively.
The set of positive reals is denoted by $\R_{+}$.
We denote by $(x_{1}, x_{2}, \dots, x_{N}) \in \R^{n_1 + n_2 + \dots + n_N}$ the column vector formed by the concatenation of vectors $x_{1} \in \R^{n_1}$ to $x_{N} \in \R^{n_N}$.
Given scalars/matrices $M_1, \dots, M_N$ (not necessarily of the same dimensions), we denote by $\texttt{diag}(M_1, \dots, M_N)$ the block diagonal matrix formed by their concatenation.
Given $A \inR{n \times m}$, $A_{i,j}$ is its $(i,  j)$-th element and $A\T$ its transposed.
Given $A \inR{n \times n}$, $A^{-1}$ is its inverse and $A^{1/2}$ is the matrix that satisfies $A = A^{1/2} A^{1/2}$ (assuming they exist).
Given $x\inR{n}$ and $A\in\Sp{n}$, we define $\|x\| \doteq \sqrt{x\T x}$, $\|x\|_A \doteq \sqrt{x\T A x}$, and $\| x \|_\infty \doteq \max_{j = 1 \dots n}{| x_{(j)} |}$, where $x_{(j)}$ is the $j$-th element of $x$.
We denote an ellipsoid defined by a given $P \in \Sp{n}$, $c \inR{n}$ and $r \in \R_{+}$ by $\E(P, c, r) \doteq \set{ x \inR{n}} {(x - c)\T P (x - c) \leq r^2 }$.
We use the shorthand $\E \equiv \E(P, c, r)$ if the values of $P$, $c$ and $r$ are clear from the context.
Given two integers $i$ and $j$ with ${j \geq i}$, $\Z_i^j$ is the set of integer numbers from $i$ to $j$.
The indicator function of a set $\cc{C} \inR{n}$ is denoted by $\delta_\cc{C} \colon \cc{C} \to \{ 0,+\infty \}$, i.e., $\delta_\cc{C}(x) = 0$ if $x \in \cc{C}$ and $\delta_\cc{C}(x) = +\infty$ if $x \not\in \cc{C}$.
The vector of ones of dimension $n$ is denoted by $\mathds{1}_n$.

\section{Introduction} \label{sec:intro}

Model Predictive Control (MPC) formulations typically rely on a terminal constraint to guarantee stability of the closed-loop system \cite{Camacho_S_2013}, \cite{Mayne_AUT_2000}, \cite{Kerrigan_Thesys_2001}, where the terminal set is taken as an Admissible Invariant Set (AIS) of the system \cite{Mayne_AUT_2000}, \cite{Kerrigan_Thesys_2001}, i.e., an invariant set that satisfies the system constraints.
The use of a terminal AIS is employed in both nominal MPC \cite[\S 2]{Rawlings_MPC_2017}, as well as in many robust MPC formulations, in which case a \textit{robust} AIS is used \cite{Chisci_AUT_01}.

Typically, in the case of linear MPC, a terminal AIS in the form of a polyhedral set is used, i.e., a set of the~form ${\{x \inR{n}: A_t x \leq b_t \}}$, where $A_t \inR{n_t \times n}$ and ${b_t \inR{n_t}}$.
In this case, the resulting optimization problem is a Quadratic Programming (QP) problem.
This is also the case when the maximal AIS is used, since it is a polyhedral set for controllable linear systems subject to linear constraints.
The computation of a polyhedral AIS for MPC is a well researched field \cite{Fiacchini_AUT_2010}, but its use typically results in the addition of a large amount of inequality constraints, i.e., $n_t$ is very large \cite[\S 5]{Blanchini_A_1999}, even for average-sized systems.
Thus, even though the resulting problem is a QP, for which many efficient solvers are available (e.g., \cite{Stellato_OSQP}, \cite{Ferreau_2014_qpOASES}), it may be resource-intensive and computationally demanding to solve.

To avoid this issue, it may be desirable to reduce the complexity of this terminal constraint. This becomes a necessity if a polyhedral AIS of the system is computationally intractable, as is often the case in many non-trivial systems, or in the robust MPC scenario.
One possible approach is to substitute the polyhedral AIS by one in the form of an ellipsoid $\E(P, c, r)$ \cite[\S 4.1]{Blanchini_A_1999}, \cite{Alvarado_JPC_RMPC_20}, \cite{Wan_AUT_03}.
This is often computationally affordable, since the ellipsoid can be obtained from the solution of a convex optimization problem subject to Linear Matrix Inequalities (LMI) whose complexity scales more favorably with the dimension of the system \cite{Boyd_SIAM_1994_LMI}, \cite{Kothare_AUT_96}.

However, the use of a terminal ellipsoidal AIS comes at the expense of the resulting MPC optimization problem no longer being a QP, due to the inclusion of a terminal quadratic constraint.
Instead, the resulting problem is a Quadratically-Constrained Quadratic Programming (QCQP) problem, which is generally more computationally demanding to solve.
Alternatively, the terminal quadratic constraint of the MPC problem can instead be posed as a Second Order Cone (SOC) constraint.
Indeed, simple algebra shows that the constraint $x \in \E(P, c, r)$, with $x \in \R^n$, can be posed as $\| P^{1/2} x - b \| \leq r$, where $b$ is the vector that solves $P^{1/2} b = P c$.
Therefore, MPC with a terminal quadratic constraint can be solved using SOC programming solvers, such as \cite{ODonoghue_SCS_21} or \cite{Garstka_JOTA_2021}.

A popular approach for solving linear MPC optimization problems is the use of operator-splitting methods, due to their simplicity and good performance.
Some solvers based on these methods are \cite{ODonoghue_SCS_21}, which uses the Douglas-Rachford method, and \cite{Stellato_OSQP, Garstka_JOTA_2021}, which use the Alternating Direction Method of Multipliers (ADMM) \cite{Boyd_FTML_2011}.
A property that they share is that they can easily consider constraints onto any closed convex set for which there is an explicit or easily implementable projection operator.
See, for instance, the approach from \cite[\S 2]{Garstka_JOTA_2021}, where the authors only consider constraints on non-empty closed convex cones, but that can also be applied to constraints on sets with a known projection operator.
Therefore, these solvers can be directly used (or easily adapted) to solve the MPC problem with a terminal quadratic constraint by posing it as the aforementioned SOC constraint.
To do so, they require $n+1$ slack variables to handle the projection onto the SOC \cite[Theorem 3.3.6]{Bauschke_Thesis_1996}.

The question that we address in this paper is if a more direct approach for considering this terminal quadratic constraint is possible, i.e., one that does not require its reformulation as a SOC constraint.
The direct inclusion of the constraint following the approach used in the aforementioned solvers is not a viable option, since it would require the projection onto the ellipsoid at each iteration of the algorithm.
The problem with this approach is that there is no explicit solution for the projection onto an ellipsoid; and would thus require the use of an iterative algorithm such as the one presented in \cite[\S 2]{Dai_SJO_2006}.
Instead, in this paper we propose an approach for the direct inclusion of ellipsoidal constraints in ADMM, without requiring additional slack variables nor its reformulation as a SOC constraint.
The main contribution is to show how a simple linear transformation of a particular part of the equality constraints of the ADMM optimization problem leads to an explicit solution of the projection step related to the ellipsoidal constraints.
The resulting projection is a generalization of the projection onto the unit ball.
The benefit of this approach is that it retains the simple matrix structures that are exploited by the solver proposed in \cite{Krupa_TCST_20}.
That is, the approach from \cite{Krupa_TCST_20} can be used to deal with the remaining equality and inequality constraints, leading to a solver with a similar small memory footprint and iteration complexity \cite{Krupa_ECC_18}.
See \cite[Chapter 5]{Krupa_Thesis_21} for a more in-depth explanation of the approach from \cite{Krupa_TCST_20}, which was also used in \cite{Krupa_TCST_21} for the development of an efficient solver for the MPC \textit{for tracking} formulation~\cite{Ferramosca_A_2009}.
The contribution of the article, when compared with \cite{Krupa_TCST_20} and \cite{Krupa_TCST_21}, is in how the ADMM algorithm is posed to deal with the terminal ellipsoidal constraint (instead of the terminal equality constraints used in \cite{Krupa_TCST_20} and \cite{Krupa_TCST_21}).
We show numerical results suggesting that the proposed approach may be preferable to transforming the resulting quadratic constraint into a SOC constraint.
The proposed solver is available in the SPCIES toolbox \cite{SPCIES} at \url{https://github.com/GepocUS/Spcies}.

The remainder of this article is structured as follows.
Section~\ref{sec:MPC} introduces the MPC formulation.
Section~\ref{sec:solver} presents the sparse ADMM-based solver for the MPC formulation presented in Section~\ref{sec:MPC}.
Section~\ref{sec:case:study} includes two case studies: one comparing the proposed solver against other alternatives from the literature, and another where we implement the solver in a low-resource embedded system to control a 12-state, 6-input chemical plant.
We close the article with Section~\ref{sec:conclusions}.

\section{Model predictive control formulation} \label{sec:MPC}

We consider the following linear MPC formulation subject to a terminal quadratic constraint:
\begin{subequations} \label{eq:MPC} 
\begin{align}
    \min\limits_{\vv{x}, \vv{u}} \;& \Sum{i = 0}{N-1} \left( \| x_i - x_r \|^2_Q + \| u_i - u_r \|^2_R \right) + \| x_{N} - x_r \|^2_T \label{eq:MPC:cost} \\
    {\rm s.t.}& \; x_0 = x(t) \label{eq:MPC:initial} \\
        & \; x_{i+1} = A x_i + B u_i, \; i\in\Z_0^{N-1} \label{eq:MPC:prediction} \\
        & \; \underline{x}_i \leq x_i \leq \overline{x}_i, \; i\in\Z_1^{N-1} \label{eq:MPC:constraint:x} \\
        & \; \underline{u}_i \leq u_i \leq \overline{u}_i, \; i\in\Z_0^{N-1} \label{eq:MPC:constraint:u}\\
        & \; x_N \in \E(P, c, r), \label{eq:MPC:terminal}
\end{align}
\end{subequations}
where $x(t) \inR{n}$ is the state, at the current discrete time instant $t$, of the system described by the controllable linear time-invariant discrete-time state-space model
\begin{equation} \label{eq:model}
    x(t+1) = A x(t) + B u(t),
\end{equation}
with $A \inR{n \times n}$, $B \inR{n \times m}$, and $u(t) \inR{m}$ is the system input applied at $t$; $x_r \inR{n}$ and $u_r \inR{m}$ are the given state and input reference, respectively, which we assume correspond to a steady state of \eqref{eq:model}; $N$ is the prediction horizon; $\vv{x} = (x_0, \dots, x_N)$ and ${\vv{u} = (u_0, \dots, u_{N-1})}$ are the predicted states and control actions throughout the prediction horizon; $Q \in \Ssp{n}$, $R \in \Ssp{m}$ and $T \in \Ssp{n}$ are the cost function matrices; $\underline{x}_i$, $\overline{x}_i \inR{n}$ and $\underline{u}_i$, $\overline{u}_i \inR{m}$ are the upper and lower bounds on the state and control input, respectively, which we assume satisfy $\underline{x}_i < \overline{x}_i$, $\underline{u}_i < \overline{u}_i$, and may be different for each prediction step $i$; and $\E(P, c, r)$ is the ellipsoid defined by the given $P \in \Sp{n}$, $c \inR{n}$ (typically $c = x_r$) and $r \in \R_{+}$.
Our consideration of step-dependent constraints in \eqref{eq:MPC:constraint:x} and \eqref{eq:MPC:constraint:u} allows for the implementation of the tightened constraints used in some tube-based robust MPC approaches \cite{Alvarado_JPC_RMPC_20}, \cite{Limon_JPC_2010}, \cite{Mayne_A_2005}.

\begin{remark}
The results presented in this paper are particularized to the MPC formulation \eqref{eq:MPC} due to its practical usefulness and because it allows us to directly apply the results of \cite{Krupa_TCST_20}.
However, they could be extended to other similar optimization problems involving quadratic constraints.
\end{remark}

\section{Sparse ADMM-based solver} \label{sec:solver}

This section presents the sparse ADMM solver for problem~\eqref{eq:MPC} in which the terminal quadratic constraint is dealt with without reformulating it into a SOC constraint.
We start by briefly recalling the ADMM algorithm, and then show how \eqref{eq:MPC} can be solved making use of the sparse approach from~\cite{Krupa_TCST_20}.

\subsection{Alternating Direction Method of Multipliers} \label{sec:ADMM}

Let $f \colon \R^{n_z} {\rightarrow} \R \cup \{+\infty\}$ and $g \colon \R^{n_v} {\rightarrow} \R \cup \{+\infty\}$ be proper, closed and convex functions, $z \inR{n_z}$, $v \inR{n_v}$, $C \inR{n_\lambda \times n_z}$ and $D \inR{n_\lambda \times n_v}$.
Consider the optimization problem%
\begin{subequations} \label{eq:optimization:problem}
    \begin{align}
        \min\limits_{z, v} &\; f(z) + g(v) \\
        {\rm s.t.} &\; C z + D v = 0, \label{eq:optimization:problem:equality}
    \end{align}
\end{subequations}
with augmented Lagrangian $\cc{L}_\rho : \R^{n_z} {\times} \R^{n_v} {\times} \R^{n_\lambda} {\rightarrow} \R$,
\begin{equation} \label{eq:lagrangian}
    \cc{L}_\rho(z, v, \lambda) {=} f(z) {+} g(v) {+} \lambda\T (Cz {+} Dv) {+} \frac{\rho}{2} \| C z {+} D v \|^2,
\end{equation}
where $\lambda \inR{n_\lambda}$ are the dual variables, and the scalar $\rho \in \R_{+}$ is the penalty parameter.
We denote a solution point of \eqref{eq:optimization:problem} by $(z^*, v^*, \lambda^*)$, assuming that one exists.

Algorithm \ref{alg:ADMM} shows the ADMM algorithm \cite{Boyd_FTML_2011} applied to problem \eqref{eq:optimization:problem} for the given exit tolerances $\epsilon_p, \epsilon_d \in \R_{+}$ and initial point $(v^0, \lambda^0)$, where the superscript $k$ denotes the value of the variable at iteration $k$ of the algorithm.
The exit condition of the algorithm is determined by the primal ($r_p$) and dual ($r_d$) residuals \cite[\S 3.3]{Boyd_FTML_2011}. These are given by
\begin{equation*}
    r_p = \| C z^{k} + D v^{k} \|_\infty,\;\; r_d = \| v^{k} - v^{k-1} \|_\infty.
\end{equation*}
The algorithm returns a suboptimal solution $(\tilde z^*, \tilde v^*, \tilde \lambda^*)$ of \eqref{eq:optimization:problem}, whose suboptimality is determined by the values of $\epsilon_p$ and $\epsilon_d$.

\begin{algorithm}[t]
    \DontPrintSemicolon
    \caption{ADMM} \label{alg:ADMM}
    \Input{$v^0$, $\lambda^0$, $\rho > 0$, $\epsilon_p > 0$, $\epsilon_d > 0$}
    $k \gets 0$\;
    \Repeat{$r_p \leq \epsilon_p$ \Kwand $r_d \leq \epsilon_d$ \label{alg:ADMM:step:exit}}{
        $z^{k+1} \gets \arg\min\limits_{z} \cc{L}_\rho(z, v^k, \lambda^k) $\; \label{alg:ADMM:step:z}
        $v^{k+1} \gets \arg\min\limits_{v} \cc{L}_\rho(z^{k+1}, v, \lambda^k) $\; \label{alg:ADMM:step:v}
        $\lambda^{k+1} \gets \lambda^k + \rho ( C z^{k+1} + D v^{k+1})$\; \label{alg:ADMM:step:lambda}
        $k \gets k + 1$\;
    }
    \KwOut{ $\tilde z^* \gets z^{k}$, $\tilde v^* \gets v^{k}$, $\tilde \lambda^* \gets \lambda^{k}$}
\end{algorithm}

\subsection{Solving the quadratically-constrained MPC using ADMM}

We cast problem \eqref{eq:MPC} as an optimization problem \eqref{eq:optimization:problem} as follows.
Let us define the auxiliary variables $\tx_i \in \R^n$, $i \in \Z_1^N$, and $\tu_i \in \R^m$, $i \in \Z_0^{N-1}$, and take
\begin{align*}
    z &= (u_0, x_1, u_1, x_2, u_2, \dots, x_{N-1}, u_{N-1}, x_{N}), \\
    v &= (\tu_0, \tx_1, \tu_1, \tx_2, \tu_2, \dots, \tx_{N-1}, \tu_{N-1}, \tx_{N}).
\end{align*}
To facilitate readability, we divide $z$ and $v$ into two parts, given by $z = (\zi, \zN)$ and $v = (\vi, \vN)$, where
\begin{align*}
    \zi &\doteq (u_0, x_1, u_1, x_2, u_2, \dots, x_{N-1}, u_{N-1}), \\
    \vi &\doteq (\tu_0, \tx_1, \tu_1, \tx_2, \tu_2, \dots, \tx_{N-1}, \tu_{N-1}),
\end{align*}
$\zN \doteq x_N$ and $\vN \doteq \tx_N$.
Then, problem \eqref{eq:MPC} can be recast as \eqref{eq:optimization:problem} by taking
\begin{subequations} \label{eq:def_fg}
\begin{align}
    f(z) &= \frac{1}{2} z\T H z + q\T z + \delta_{(G z - b = 0)}(z), \label{eq:def_fg:f} \\
    g(v) &= \delta_{(\lbvi \leq \vi \leq \ubvi)}(\vi) + \delta_{\E(P, c, r)}(\vN), \label{eq:def_fg:g}
\end{align}
\end{subequations}
where $H$, $q$ account for the cost function \eqref{eq:MPC:cost}, $\lbvi$, $\ubvi$ for the inequality constraints \eqref{eq:MPC:constraint:x} and \eqref{eq:MPC:constraint:u}, and $G$, $b$ for the equality constraints \eqref{eq:MPC:initial} and \eqref{eq:MPC:prediction}.
The common and straightforward approach is to take the ADMM equality constraint \eqref{eq:optimization:problem:equality} as $z - v = 0$.
However, to handle the ellipsoidal constraint we propose to take instead
\begin{subequations} \label{eq:constraints:ADMM}
    \begin{align}
        &\zi - \vi = 0, \label{eq:constraints:ADMM:i} \\
        &P^{1/2} (\zN - \vN) = 0, \label{eq:constraints:ADMM:N}
    \end{align}
\end{subequations}
where $P^{1/2}$ is the matrix that satisfies $P = P^{1/2} P^{1/2}$ and our choice of imposing \eqref{eq:constraints:ADMM:N} instead of the more straightforward $\zN - \vN = 0$ is the key for being able to efficiently deal with the quadratic constraint.
That is, we take matrices $C$ and $D$~as
\begin{align*}
    C &= \texttt{diag}(I_m,\, I_n,\, I_m,\, \dots, I_n,\, I_m,\, P^{1/2}), \\
    D &= -\texttt{diag}(I_m,\, I_n,\, I_m,\, \dots, I_n,\, I_m,\, P^{1/2}).
\end{align*}
Let us also define $\lambda = (\lambdai, \lambdaN)$, where $\lambdai$ are the dual variables associated to the constraints \eqref{eq:constraints:ADMM:i}, and $\lambdaN$ are the dual variables associated to the constraints \eqref{eq:constraints:ADMM:N}.

\begin{algorithm}[t]
    \DontPrintSemicolon
    \caption{Sparse ADMM-based solver for \eqref{eq:MPC}} \label{alg:ADMM:MPC}
    \Input{$x(t)$, $(x_r, u_r)$, $v^0$, $\lambda^0$, $\epsilon_p > 0$, $\epsilon_d > 0$}
    $q \gets -(R u_r, Q x_r, R u_r, \dots, Q x_r, R u_r, T x_r)$ \label{alg:ADMM:MPC:step:q} \;
    $b \gets (-A x(t), 0, 0, \dots, 0)$, $k \gets 0$ \label{alg:ADMM:MPC:step:b} \;
    \Repeat{$r_p \leq \epsilon_p$ \Kwand $r_d \leq \epsilon_d$ \label{alg:ADMM:MPC:step:exit}}{
    $\hat q_k \gets q + (\lambdai^k - \rho \vi^k, P^{1/2} \lambdaN^k - \rho P \vN^k )$\label{alg:ADMM:MPC:step:qhat}\;
    $y \gets -(G \hat H^{-1} \hat q_k + b)$\; \label{alg:ADMM:MPC:step:rhs}
    $\mu \gets$ solution of $G \hat{H}^{-1} G\T \mu {=} y$ using \cite[Alg. 11]{Krupa_Thesis_21}\; \label{alg:ADMM:MPC:step:mu}
    $z^{k+1} \gets - \hat H^{-1} ( G\T \mu + \hat q_k)$\; \label{alg:ADMM:MPC:step:z}
    $\vi^{k+1} \gets \max \{ \min \{ \zi^{k+1} + \rho^{-1} \lambdai^k, \; \ubvi \}, \; \lbvi \}$ \label{alg:ADMM:MPC:step:vi}\;
    $\vN^{k+1} \gets \zN^{k+1} + \rho^{-1} P^{-1/2} \lambdaN^k$ \label{alg:ADMM:MPC:step:vN}\;
    \If{$(\vN^{k+1} - c)\T P (\vN^{k+1} - c) > r^2$}{
        $\vN^{k+1} \gets \fracg{r (\vN^{k+1} - c)}{\sqrt{(\vN^{k+1} - c)\T P (\vN^{k+1} - c)}} + c$\;
    } \label{alg:ADMM:MPC:step:vN:end}
    $\lambdai^{k+1} \gets \lambdai^k + \rho ( \zi^{k+1} - \vi^{k+1})$ \label{alg:ADMM:MPC:step:lambdai}\;
    $\lambdaN^{k+1} \gets \lambdaN^k + \rho P^{1/2} ( \zN^{k+1} - \vN^{k+1})$ \label{alg:ADMM:MPC:step:lambdaN}\;
    $k \gets k + 1$\;
    }
    \KwOut{$\tilde z^* \gets z^{k}$, $\tilde v^* \gets v^{k}$, $\tilde \lambda^* \gets \lambda^{k}$}
\end{algorithm}

By taking the ingredients of problem \eqref{eq:optimization:problem} this way, the steps of Algorithm~\ref{alg:ADMM} can be computed efficiently.
This leads to Algorithm~\ref{alg:ADMM:MPC}, whose steps we explain in detail in the following.

From \eqref{eq:def_fg:f}, we have that Step \ref{alg:ADMM:step:z} of Algorithm~\ref{alg:ADMM} requires solving the equality-constrained QP problem
\begin{subequations} \label{eq:QP:f}
\begin{align}
    \min\limits_{z} \; &\frac{1}{2}\, z\T \hat{H} z + \hat{q}_k\T z \\
    {\rm s.t.} \; & G z = b,
\end{align}
\end{subequations}
where
\begin{align*}
    \hat H &= H + \rho C\T C = H + \rho \cdot \texttt{diag}(I_m, I_n, \dots, I_n, I_m, P), \\
    \hat q_k &= q + (\lambdai^k - \rho \vi^k, P^{1/2} \lambdaN^k - \rho P \vN^k ),
\end{align*}
are obtained by adding to \eqref{eq:def_fg:f} the terms in the Lagrangian~\eqref{eq:lagrangian}.
There are several different ways in which \eqref{eq:QP:f} could be explicitly solved.
However, the matrices of \eqref{eq:QP:f} have a particular sparse banded structure that commonly arises in linear MPC.
Thus, we consider the use of the procedure presented in \cite[\S IV.A]{Krupa_TCST_20}, where the authors proposed an efficient solver for this class of QP problems that exploits the banded structure of its ingredients.
Due to space considerations, we do not repeat the complete procedure here. 
Instead, the reader is referred to \cite[\S IV.A]{Krupa_TCST_20} and \cite[\S 5.1.2]{Krupa_Thesis_21} for an in-depth explanation of the approach, whose main benefit is its small computational and memory requirements \cite{Krupa_TCST_20, Krupa_ECC_18}.
Steps~\ref{alg:ADMM:MPC:step:q} and \ref{alg:ADMM:MPC:step:b} of Algorithm~\ref{alg:ADMM:MPC} compute the ingredients $q$ and $b$ of \eqref{eq:def_fg:f} for the current system state $x(t)$ and reference $(x_r, u_r)$.
Steps~\ref{alg:ADMM:MPC:step:qhat}--\ref{alg:ADMM:MPC:step:z} of Algorithm~\ref{alg:ADMM:MPC} perform the procedure from \cite[\S IV.A]{Krupa_TCST_20}.
Step~\ref{alg:ADMM:MPC:step:z} solves the linear system $G \hat{H}^{-1} G\T \mu = y$ using \cite[Alg. 11]{Krupa_Thesis_21} by means of the Cholesky decomposition of $G \hat{H}^{-1} G\T$, which is computed offline.
All matrix operations in Steps~\ref{alg:ADMM:MPC:step:rhs}--\ref{alg:ADMM:MPC:step:z} of Algorithm~\ref{alg:ADMM:MPC} are performed sparsely by directly exploiting the matrix structures.

Step \ref{alg:ADMM:step:v} of Algorithm~\ref{alg:ADMM} has a separable structure that allows the problem to be divided into two parts: one for $\vi$ and one for $\vN$.
The update of $\vi^{k+1}$ is the solution of the QP problem
\begin{subequations} \label{eq:QP:g0}
\begin{align}
    \min\limits_{\vi} \; &\frac{1}{2}\, \vi\T \vi - (\zi^{k+1} + \rho^{-1} \lambdai^k)\T \vi \\
    {\rm s.t.} \; &\lbvi \leq \vi \leq \ubvi.
\end{align}
\end{subequations}
Note that each element of $\vi$ in \eqref{eq:QP:g0} is separable, constrained to a box and has a quadratic cost.
Thus, \eqref{eq:QP:g0} has the well-known explicit solution provided in Step~\ref{alg:ADMM:MPC:step:vi} of Algorithm~\ref{alg:ADMM:MPC}~\cite[\S III]{Krupa_TCST_20}.
The update of $\vN^{k+1}$ is the solution of
\begin{align*}
    \min\limits_{\vN}\;& \frac{\rho}{2} \vN\T P \vN - ( \rho P \zN^{k+1} + P^{1/2} \lambdaN^k )\T \vN \\
    {\rm s.t.} \;& \vN \in \E(P, c, r),
\end{align*}
which, dividing the objective function by $\rho$ and denoting $P^{-1/2} = P^{-1} P^{1/2}$, can be recast as
\begin{subequations} \label{eq:QP:vN}
    \begin{align}
        \min\limits_{\vN}\;& \frac{1}{2} \| \vN {-} ( \zN^{k+1} + \rho^{-1} P^{-1/2} \lambdaN^k ) \|^2_P \\
        {\rm s.t.} \;& \vN \in \E(P, c, r),
    \end{align}
\end{subequations}
whose explicit solution is given by the following proposition.

\begin{proposition} \label{theo:P:projection}

Let $a \inR{n}$, $P \in\Sp{n}$, $c \inR{n}$, $r \in \R > 0$. Then, the solution $v^*$ of the convex optimization problem
\begin{equation}
    \min\limits_{v}\; \left\{ \frac{1}{2} \| v - a \|^2_P, \; {\rm s.t.} \; v \in \E(P, c, r) \right\}, \label{eq:P:projection:optimization:problem}
\end{equation}
is given by
$v^* = \begin{cases} \quad\quad\quad\quad a &\text{if~} a \in \E(P, c, r) \vspace*{0.5em}\\
\fracg{r (a {-} c)}{\sqrt{ (a {-} c)\T P (a {-} c)}} + c &\text{otherwise.} \end{cases}$
\end{proposition}

\begin{remark}
Proposition~\ref{theo:P:projection} is a generalization of the well-known projection onto the unit ball \cite[\S 6.5.1]{Parikh_proximal_2014}, for values of $P$ other than the identity matrix.
In fact, the proposition can be proven by taking the change of variable $\tilde{v} = r^{-1} P^{1/2} (v - c)$, which is well-defined since $r > 0$.
In Appendix~\ref{app:proof:theo:P:projection} we provide an alternative proof based on analyzing the dual problem of~\eqref{eq:P:projection:optimization:problem}.
\end{remark}

Steps~\ref{alg:ADMM:MPC:step:vN}--\ref{alg:ADMM:MPC:step:vN:end} of Algorithm~\ref{alg:ADMM:MPC} solve \eqref{eq:QP:vN} using Proposition~\ref{theo:P:projection}.
Finally, Steps~\ref{alg:ADMM:MPC:step:lambdai} and \ref{alg:ADMM:MPC:step:lambdaN} of Algorithm~\ref{alg:ADMM:MPC} perform Step~\ref{alg:ADMM:step:lambda} of Algorithm~\ref{alg:ADMM}.
The control action $u(t)$ to be applied to the system is taken as the first $m$ elements of $\tilde v^*$.

\begin{remark} \label{rem:P:projection} 
Our imposition of \eqref{eq:constraints:ADMM:N}, instead of the more immediate $\zN - \vN = 0$, is the key that allows us to efficiently deal with the quadratic constraint, since it leads to the optimization problem with explicit solution \eqref{eq:QP:vN}.
The use of the typical constraint $\zN - \vN = 0$ would have instead resulted in \eqref{eq:QP:vN} being an Euclidean projection onto the ellipsoid, which, as mentioned in the introduction, does not have an explicit solution.
However, by multiplying the ADMM equality constraints related to the quadratic constraint by $P^{1/2}$, we can instead use Proposition~\ref{theo:P:projection}.
Note that pre-multiplying the constraint by $P^{1/2}$ does not affect the optimal solution of the optimization problem, since $\zN - \vN = 0$ is still the only solution of \eqref{eq:constraints:ADMM:N}.
\end{remark}

\section{Case study} \label{sec:case:study}

We present two case studies: one comparing the proposed solver against other state-of-the-art solvers, and another where we implement it on an embedded system to control a 12-state, 6-input reactor plant.

\subsection{Comparison with state-of-the-art solvers} \label{sec:case:study:comparison}

\begin{figure}[t]
    \centering
    \includegraphics[width=0.9\columnwidth]{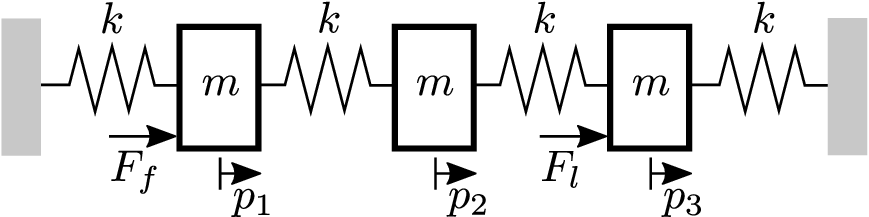}
    \caption{Chain of three masses.}
    \label{fig:masses}
\end{figure}

We first consider a system of three objects connected by springs, illustrated in Fig. \ref{fig:masses}, inspired by the case study from \cite[\S 5]{Kogel_IFAC_2011}.
We take the mass of the outer objects as ${m = 1}$kg, and the mass of the central object as $m = 0.5$kg. The spring constants are all taken as $k = 2$N/m. There are two external forces acting on the system: a force $F_f$ [N] acting on the first object, and a force $F_l$ [N] acting on the last object. The state of the system is given by the position, $p_i$ [dm]\footnote{We consider positions in decimeters to improve the numerical conditioning.}, and velocity, $v_i$ [m/s], of each object, i.e., $x = (p_1, \, p_2, \, p_3, \, v_1, \, v_2, \, v_3)$; and the input is $u = (F_f, \, F_l)$. We compute a model \eqref{eq:model} by taking a $0.2$s sampling time. We consider the following constraints:
\begin{equation} \label{eq:masses:constraints}
    -10 \leq p_i \leq 3\, \text{dm},\; i \in \Z_1^3, \; | F_f | \leq 0.8\, \text{N}, \; | F_l | \leq 0.8\, \text{N}.
\end{equation}
We design an MPC controller \eqref{eq:MPC} with the parameters $N = 10$, $Q = \texttt{diag}(15, \, 15, \, 15, \, 1, \, 1, \, 1)$ and $R = 0.1 I_2$.
We select the reference as $x_r = (2.5, 2.5, 2.5, 0, 0, 0)$ and $u_r = (0.5, 0.5)$.
We take $c = x_r$, $r = 1$, and jointly compute a matrix $P$ and a state-feedback control gain $K_t$ such that the resulting ellipsoid $\E(P, c, r)$ is an AIS of the closed-loop system using standard procedures \cite{Chen_ACC_2001}, \cite[\S C.8.1]{Tarbouriech_S_2011}, \cite{Boyd_SIAM_1994_LMI}.
In particular, we follow the procedure described in \cite[Appendix B]{Krupa_arXiv_ellipMPC_21_v2}, where we find a feasible solution of the LMI optimization problem using the YALMIP toolbox \cite{YALMIP} along with the SDPT3 solver \cite{SDPT3_99}, taking the $\lambda$ parameter from \cite[Appendix B]{Krupa_arXiv_ellipMPC_21_v2} as $\lambda = 0.999$.
Finally, we take $T$ as the solution of the Lyapunov equation
\begin{equation} \label{eq:T:Lyapunov}
    (A + B K_t)\T T (A + B K_t) - T = - Q - K_t\T R K_t.
\end{equation}

{\renewcommand{\arraystretch}{1.00}%
    \begin{table*}[t]
    \setlength{\tabcolsep}{4pt}
    \centering
	\begin{threeparttable}
    \begin{tabular}{rccccccccccc}
        & \multicolumn{4}{c}{Number of iterations} & \multicolumn{4}{c}{Computation time [ms]} & \\
        \cmidrule(lr){2-5}\cmidrule(lr){6-9}\cmidrule(lr){10-10}\cmidrule(lr){11-11}
        & Average & Median & Maximum & Minimum & Average & Median & Maximum & Minimum & Param. $\rho$ & Formulation \\
        \cmidrule(lr){2-5}\cmidrule(lr){6-9}\cmidrule(lr){10-10}\cmidrule(lr){11-11}
        Algorithm \ref{alg:ADMM:MPC} & 1014.64 & 901.0 & 3035 & 128 & 1.07 & 0.95 & 4.07 & 0.13 & 280 & \eqref{eq:MPC} \\
        SOC-version & 909.72 & 713.0 & 4177 & 102 & 2.65 & 2.07 & 12.20 & 0.29 & 190 & \eqref{eq:MPC} \\
        SCS & 338.90 & 250.0 & 25700 & 25 & 2.45 & 1.81 & 182.63 & 0.20 & - & \eqref{eq:MPC} \\
        HPIPM & 22.30 & 22.0 & 49 & 12 & 0.49 & 0.48 & 1.33 & 0.27 & - & \eqref{eq:MPC} \\
        laxMPC & 626.04 & 614 & 4895 & 59 & 0.64 & 0.63 & 4.82 & 0.06 & 100 & \cite[Eq. (9)]{Krupa_TCST_20} \\
        OSQP & 378.95 & 250.0 & 7375 & 50 & 2.95 & 2.07 & 56.97 & 0.39 & - & Variant $^\star$ of \eqref{eq:MPC} \\
        FalcOpt & 2591.78 & 2692.0 & 5203 & 121 & 4.73 & 4.89 & 11.70 & 0.24 & - & Variant $^\dagger$ of \eqref{eq:MPC} \\
        \cmidrule(lr){2-11}
    \end{tabular}
	\begin{tablenotes}[] \footnotesize
        \item $^\star$ The maximal AIS of the system is used as the terminal set, instead of an ellipsoidal AIS.
        $^\dagger$ No state constraints.
        \item When applicable, the penalty parameter $\rho$ of ADMM algorithms is selected as the one that results in the lowest maximum iterations.
        \item The results correspond to the $1451$ random states $x(t)$ for which all MPC formulations were feasible.
	\end{tablenotes}
    \caption{Comparison between solvers for $2000$ randomly generated states $x(t)$ within the system constraints.}
    \label{tab:comparison}
	\end{threeparttable}
\end{table*}}

{\renewcommand{\arraystretch}{1.0}%
    \begin{table*}[t]
    \setlength{\tabcolsep}{4pt}
    \centering
	\begin{threeparttable}
    \begin{tabular}{rcccccccccccc}
        & \multicolumn{4}{c}{Number of iterations} & \multicolumn{4}{c}{Computation time [ms]} & \\
        \cmidrule(lr){2-5}\cmidrule(lr){6-9}\cmidrule(lr){10-10}
        & Average & Median & Maximum & Minimum & Average & Median & Maximum & Minimum & Param. $\rho$ \\
        \cmidrule(lr){2-5}\cmidrule(lr){6-9}\cmidrule(lr){10-10}
        Algorithm \ref{alg:ADMM:MPC} & 170.46 & 137.0 & 1198 & 2 & 0.35 & 0.21 & 2.75 & 0.01 & 50 \\
        SOC-version                  & 182.43 & 143.0 & 1198 & 2 & 0.78 & 0.47 & 6.44 & 0.01 & 50 \\
        SCS                          & 156.96 & 125.0 & 2675 & 25 & 2.06 & 0.86 & 53.34 & 0.05 & - \\
        HPIPM                        & 6.85 & 7.0 & 13 & 4 & 0.21 & 0.18 & 0.65 & 0.04 & - \\
        \cmidrule(lr){2-10}
    \end{tabular}
	\begin{tablenotes}[] \footnotesize
    \item The table shows the computational results for the $539$ feasible optimization problems out of the $1000$ tests.
	\end{tablenotes}
    \caption{Comparison between solvers applied to \eqref{eq:MPC} for $1000$ randomly generated systems.}
    \label{tab:random}
	\end{threeparttable}
\end{table*}}

\begin{figure*}[t]
    \centering
    \begin{subfigure}[ht]{0.40\textwidth}
        \includegraphics[width=0.99\linewidth]{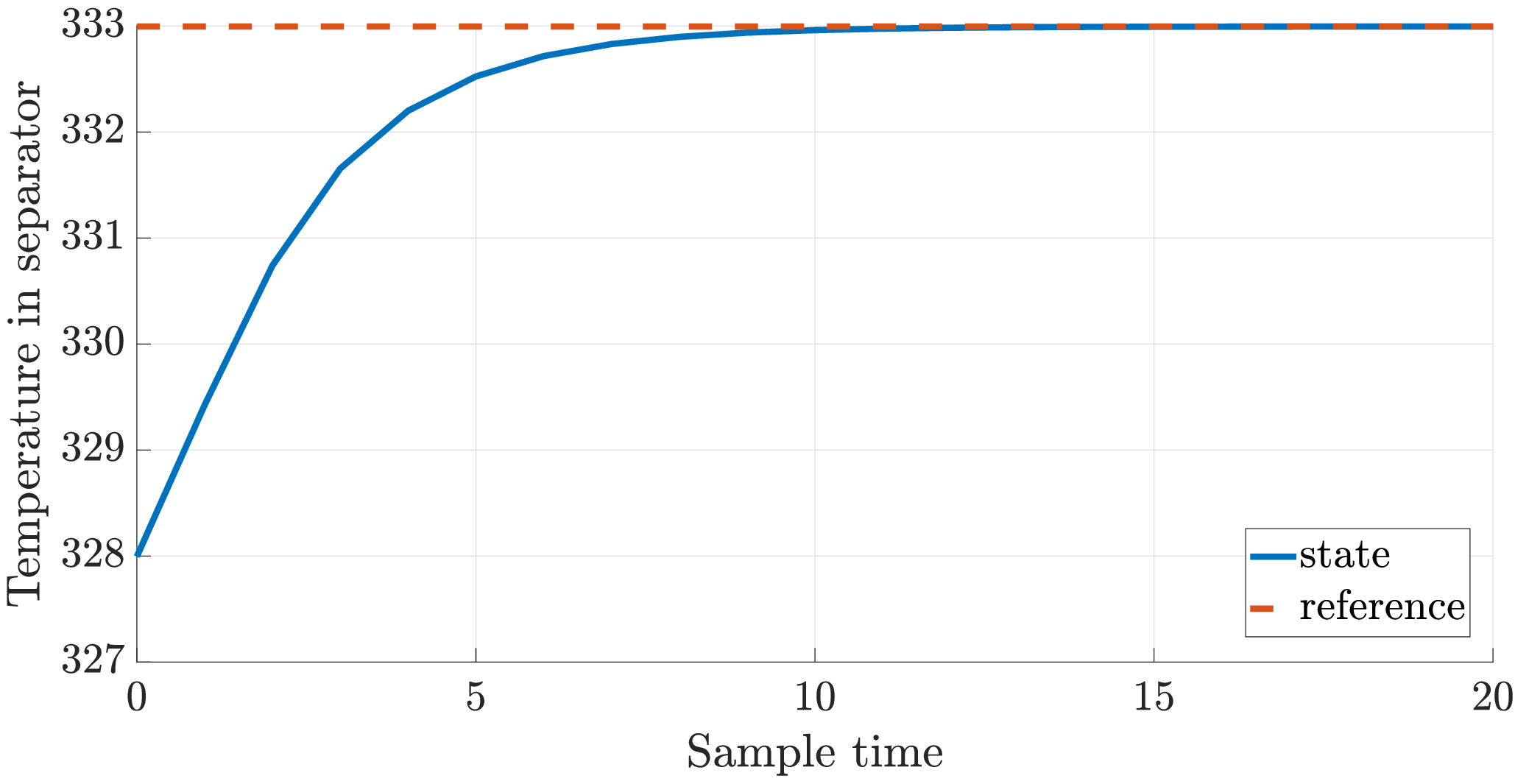}
        \caption{Temperature in the separator.}
        \label{fig:result:state}
    \end{subfigure}%
    \hspace{5em}
    \begin{subfigure}[ht]{0.40\textwidth}
        \includegraphics[width=0.99\linewidth]{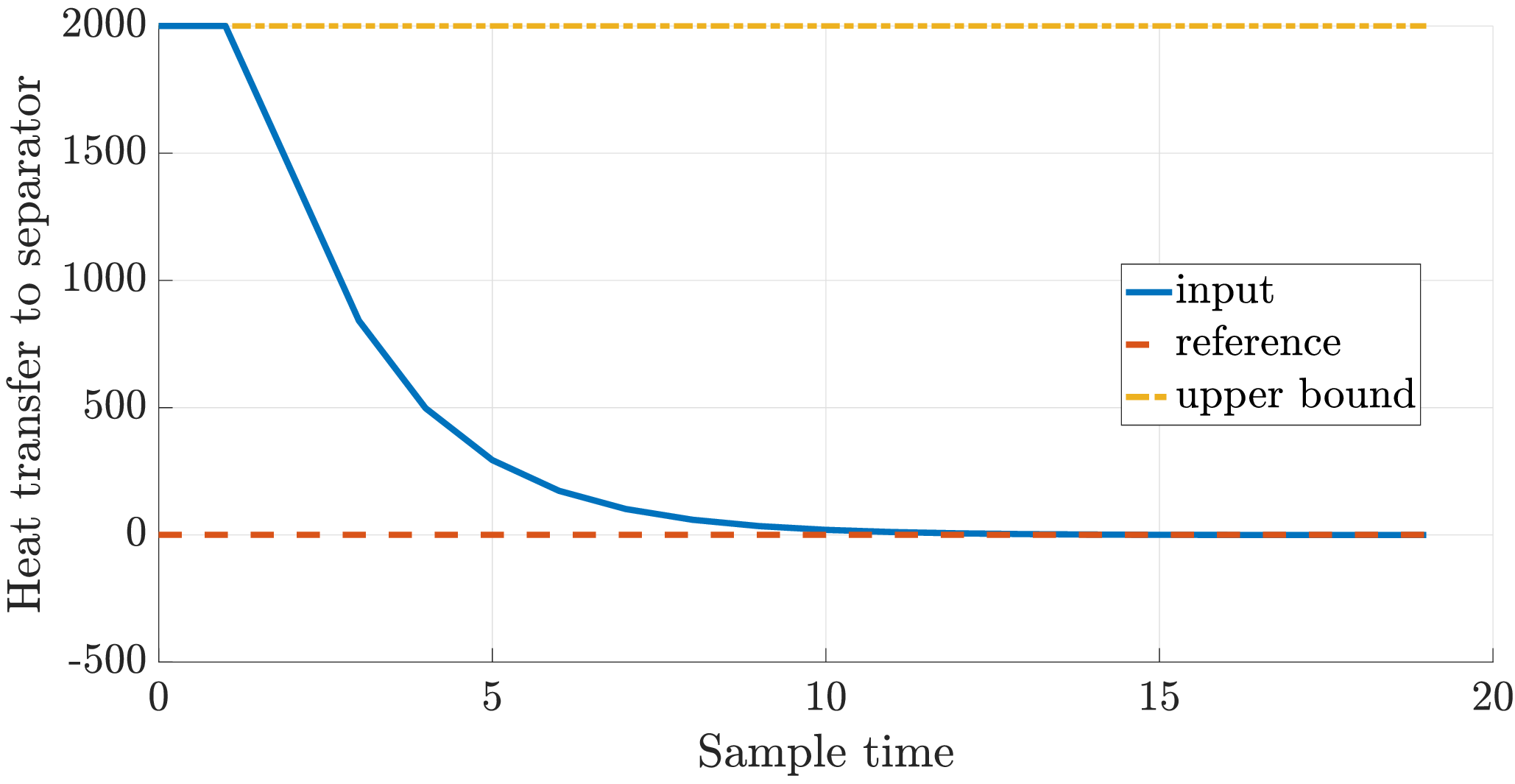}
        \caption{Heat transfer to the separator.}
        \label{fig:result:input}
    \end{subfigure}%

    \begin{subfigure}[ht]{0.40\textwidth}
        \includegraphics[width=0.99\linewidth]{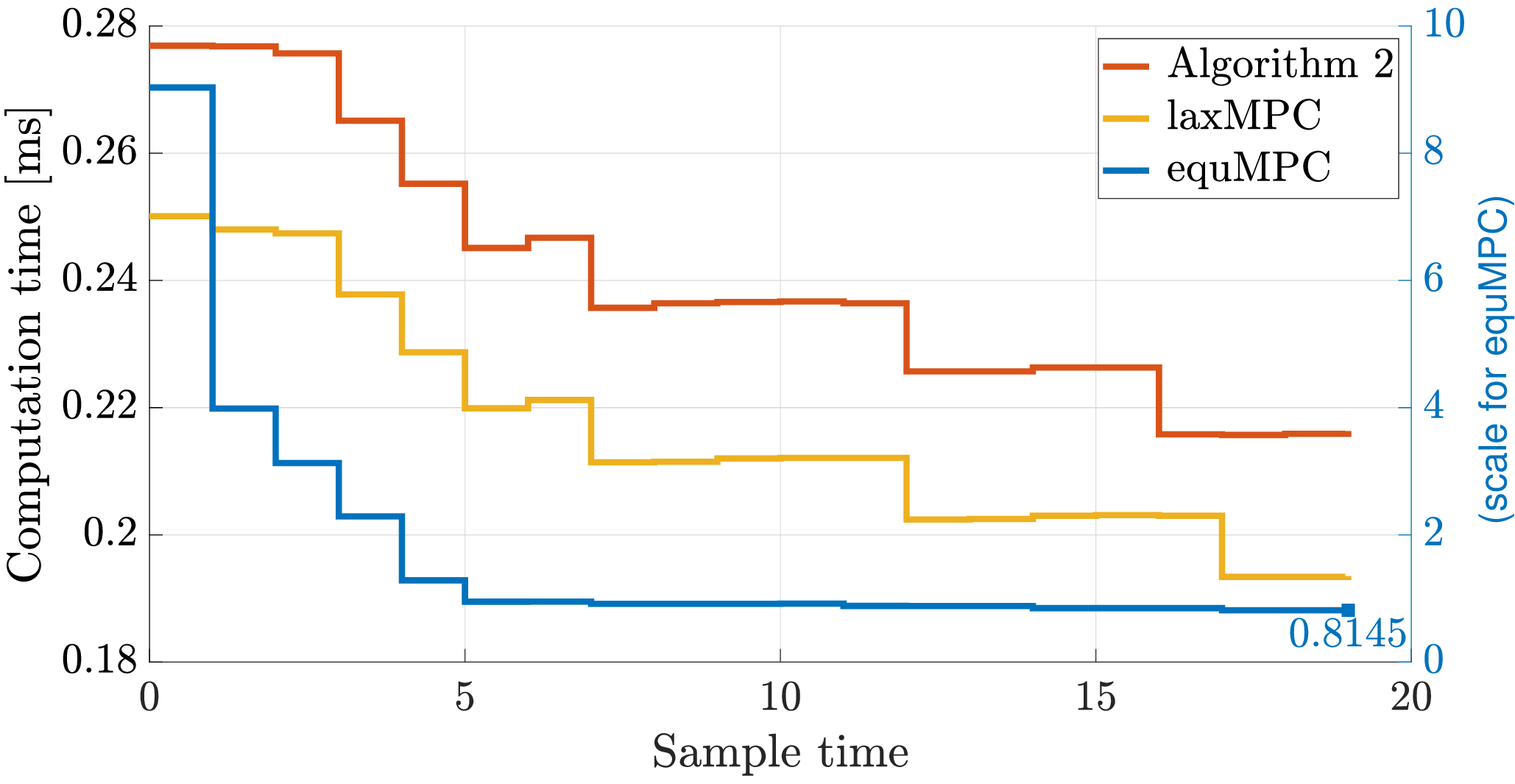}
        \caption{Computation time.}
        \label{fig:result:time}
    \end{subfigure}%
    \hspace{5em}
    \begin{subfigure}[ht]{0.40\textwidth}
        \includegraphics[width=0.99\linewidth]{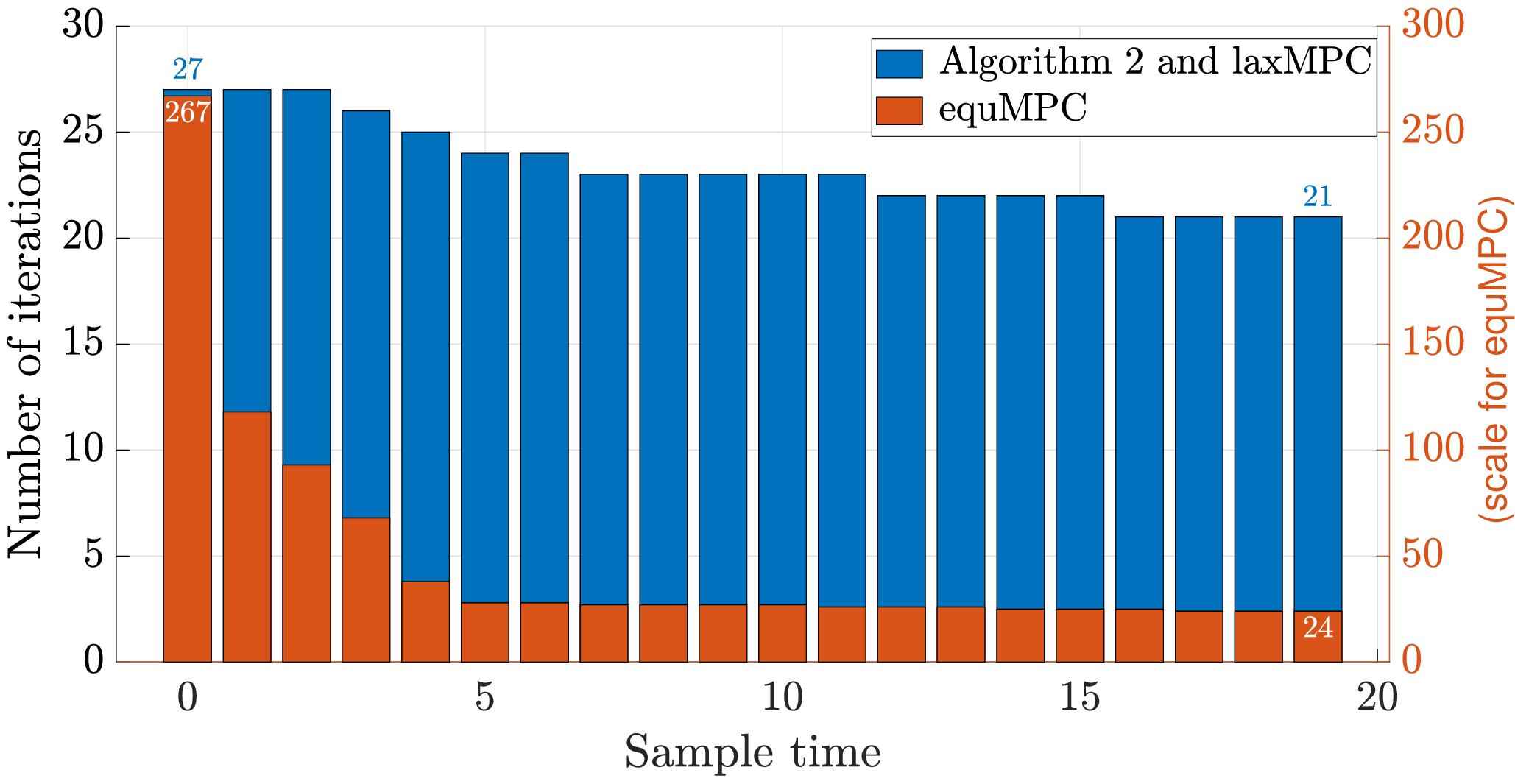}
        \caption{Number of iterations.}
        \label{fig:result:iter}
    \end{subfigure}%
    \caption{Closed-loop results on the \textit{reactors and separator} chemical system. In subfigures (c) and (d), the results for the \textit{equMPC} are presented in the right-hand-side axis.}
    \label{fig:result}
\end{figure*}

We solve problem \eqref{eq:MPC} using the proposed solver for $2000$ random states $x(t)$ taken from a uniform distribution in the~set
\begin{equation*}
    \set{x = (p_1, p_2, p_3, v_1, v_2, v_3)}{0 \leq p_i \leq 3, |v_i| \leq 0.4, i \in \Z_1^3}.
\end{equation*}
Table \ref{tab:comparison} shows the results obtained using Algorithm \ref{alg:ADMM:MPC}, taking $v^0 = 0$ and $\lambda^0 = 0$.
The table also shows the results of other alternative solvers and approaches.
In particular, we compare Algorithm \ref{alg:ADMM:MPC} with:
\begin{itemize}[leftmargin=*]
    \setlength\itemsep{0.0em}
    \item \textbf{SOC-version}. This solver transforms the quadratic constraint into a SOC constraint by following the approach taken in \cite{krupa_HMPC_solver_2023}, which resembles the approach from the COSMOS solver \cite{Garstka_JOTA_2021}.
        The resulting ADMM-based solver requires the inclusion of an additional decision variable (and corresponding dual variable) to deal with the SOC constraint.
        Additionally, in order to deal with the SOC constraint, the ADMM problem is posed in a different way to the one explained in Section~\ref{sec:solver}.
        Instead of taking $v$ as an auxiliary copy of $z$, and imposing $z = v$ in \eqref{eq:optimization:problem:equality}, the idea is to consider $v$ as ``slack" variables for the constraints of \eqref{eq:MPC}.
        The approach is very similar, but an important consequence is that the solver does not retain the simple structures exploited by the solver from \cite{Krupa_TCST_20}.
        Instead, the equality-constrained subproblem of the ADMM algorithm is solved using an LDL decomposition stored using standard sparse matrix representations.
        The results in Table~\ref{tab:comparison} show how the loss of the simple matrix structures exploited by the solver from \cite{Krupa_TCST_20} may lead to worse computation times.
        This solver is also available in the SPCIES toolbox \cite{SPCIES}.
\item \textbf{SCS} \cite{ODonoghue_SCS_21}, version \texttt{3.1.0}. A recent sparse SOC solver with quadratic cost function.
    We use it to solve problem \eqref{eq:MPC} but imposing the terminal quadratic constraint using a SOC.
    We include it to show a comparison between our proposed solver with an efficient SOC programming solver from the literature.
    The results indicate that our proposed solver is faster than SCS, mainly due to the fact that the proposed solver is particularized to our MPC formulation while SCS is for more generic optimization problems.
    \item \textbf{HPIPM} \cite{Frison_IFAC_2020}, GitHub commit \texttt{d0df588}. An interior-point QP solver ``with a particular focus on MPC" that has recently been extended to QCQP problems \cite{Frison_arXiv_2021}.
        The results using this solver outperform the other candidates, which is not surprising taking into account that its linear algebra operations are performed using the highly-efficient BLASFEO package \cite{Frison_BLASFEO_2018} (which is particularly efficient when dealing with matrices with dimensions in the few-hundreds due to its efficient cache usage).
        This solver is only included to provide a comparison with an interior point solver that is also able to handle SOC constraints (although the comparison may be unfair due to the use of BLASFEO in the HPIPM solver versus the plain-C implementation used by all other solvers).
        However, an in-depth analysis and discussion of the advantages/disadvantages of ADMM-based solvers versus interior-point solvers is beyond the scope of this paper.
    \item \textbf{laxMPC}. The sparse ADMM-based solver presented in \cite{Krupa_TCST_20} for formulation \cite[Eq. (9)]{Krupa_TCST_20}, which is exactly the same as \eqref{eq:MPC} but without a terminal constraint. This solver also uses the procedure from \cite{Krupa_TCST_20} to deal with the equality and box constraints. We included it to show the effect that the addition of the terminal quadratic constraint has on the performance of the solver.
The results show that the influence of the terminal constraint on the computation time of the solver is very small.
The solver is also available in the SPCIES toolbox \cite{SPCIES}.
    \item \textbf{OSQP} \cite{Stellato_OSQP}, version \texttt{0.6.0}. A sparse QP solver that is not particularized to MPC but that is also based on the ADMM algorithm. We use it to solve problem \eqref{eq:MPC} but with a polyhedral terminal set, thus resulting in a QP problem.
    The terminal set is taken as the maximal AIS of the system, which is computed using the MPT3 toolbox \cite{MPT3}. The resulting polyhedron $\{x \inR{n} :A_t x \leq b_t \}$ is described by a matrix $A_t$ with $274$ rows, resulting in as many inequality constraints in the QP problem.
    The results show that, even for a small-scale system, the inclusion of a polyhedral terminal AIS can significantly increase the computational time when compared to a terminal quadratic constraint.
    \item \textbf{FalcOpt} \cite{FalcOpt, FalcOpt_Theory}, GitHub commit \texttt{5ac104c}. A solver for non-linear MPC that is suitable for embedded systems and considers a terminal quadratic constraint. However, it does not consider state constraints and the use of (a variation of) the SQP method may make it less efficient when dealing with the linear MPC case than a tailored linear MCP solver.
\end{itemize}

We note that different MPC formulations are considered by the above solvers, as indicated in the column ``Formulation'' of Table~\ref{tab:comparison}.
Thus, a state $x(t)$ may lead to a feasible MPC problem for some of the solvers but not all.
Therefore, to provide a fair comparison, Table~\ref{tab:comparison} shows results for the $1451$ states for which all the considered MPC formulations were feasible.

The exit tolerances of the solvers are set to $10^{-4}$ except for the FalcOpt solver, where it is set to $10^{-1}$.
The warmstart procedure from OSQP and the acceleration and verbose options of SCS are disabled to provide a fair comparison with the other solvers.
The maximum number of iterations is set to $30000$ for all solvers.
All other options are set to their default values.
We use version \texttt{v0.3.10} of the SPCIES toolbox.
When applicable, solvers are compiled using the \texttt{gcc} compiler with the \texttt{-O3} flag.
The tests are performed using the MATLAB interface of the solvers on an Intel Core i5-8250U operating at $1.60$GHz.

We next compare the performance of the proposed approach when applied to $1000$ randomly generated systems.
We generate systems using the \texttt{drss()} function from Matlab.
The number of states $n$ and inputs of each system are selected from a uniform distribution in the intervals $[2, 10]$ and $[2, n]$, respectively.
We take the constraint as $\overline{x} = 1.5 \mathds{1}_n$, $\underline{x} = -0.3 \mathds{1}_n$, $\overline{u} = 0.3 \mathds{1}_m$, $\underline{u} = -\overline{u}$.
We set the reference to $x_r = 0$, $u_r = 0$.
For each system, we build the MPC controller~\eqref{eq:MPC} taking $N = n+3$ and $Q$ and $R$ as diagonal matrices with non-zero elements taken from a uniform distribution on the intervals $(0, 5]$ and $(0, 2]$, respectively.
As before, we construct the terminal ellipsoidal constraint using the procedure from \cite[Appendix B]{Krupa_arXiv_ellipMPC_21_v2} and take $T$ as the solution of~\eqref{eq:T:Lyapunov}.
For each system, we solve the MPC formulation~\eqref{eq:MPC} for the initial state $x(t) = \mathds{1}_n$.
We only consider the solvers capable of solving \eqref{eq:MPC}, i.e., Algorithm~\ref{alg:ADMM:MPC}, \textit{SOC-version}, SCS and HPIPM.
We take the same solver exit tolerances and options detailed in the previous paragraph.
Table~\ref{tab:random} shows the computational results for the $539$ feasible MPC problems out of the $1000$ randomly generated systems.
The value of $\rho$ for Algorithm~\ref{alg:ADMM:MPC} and \textit{SOC-version} where selected as the ones that resulted in the lowest maximum number of iterations when tested on a different set of $300$ randomly generated systems.
The results show that the proposed approach may lead to significant computational advantages when compared with the alternative of using SOC constraint.
We note that the average computational time of Algorithm~\ref{alg:ADMM:MPC} is remarkably close to HPIPM, in spite of being programmed in library-free C.

\subsection{Implementation on an embedded system} \label{sec:case:study:embedded}

We now consider the 12-state, 6-input reactor system described in \cite[\S VI.A]{Krupa_TCST_20} (see \cite[\S 5.7.1]{Krupa_Thesis_21} for a more in-depth description), in which two chemical reactions take place in two reactors and one phase separator.
The states of the system are the liquid height and temperature of each of the three volumes, as well as the concentration of the two reactants in each of them.
The inputs are the two flows of reactant fed to each reactor, the \textit{discarded flow} from the separator, and the three independent heats transferred by the cooling system to each of the volumes (the two reactors and the separator).

We implement \eqref{eq:MPC} on a $1.8$GHz Raspberry Pi 4 Model B running a 32-bit operating system using the 6.1.21 Linux kernel.
We take the following ingredients and parameters:
\begin{align*}
    Q &= \texttt{diag}(1, 1, 1, 10, 1, 1, 1, 10, 1, 1, 1, 10), \; R = 0.1 I_6, \\
    N &= 5, \; \epsilon_p = \epsilon_d = 10^{-6}, \; \rho = 0.2.
\end{align*}
We take the reference as the operating point described in \cite[Table~II]{Krupa_TCST_20}.
We compute an ellipsoidal AIS of the system centered around this reference using the same approach used in Section~\ref{sec:case:study:comparison}, but taking $T$ as the $P$ matrix of the ellipsoid.
In this case, we are unable to compute an AIS of the system using the MPT3 toolbox \cite{MPT3} due to the size of the system.

Figure~\ref{fig:result} shows the resulting trajectory of the temperature and heat transfer of the separator (12th state and 6th input, respectively), as well as the number of iterations and computation time of the proposed solver at each sample time.
Figures~\ref{fig:result:time} and \ref{fig:result:iter} also show the computation times and number of iterations when using the laxMPC formulation described in Section~\ref{sec:case:study:comparison} and the ``equMPC" formulation, which is another MPC formulation offered in the SPCIES toolbox \cite{SPCIES} that includes a terminal equality constraint.
Tests where executed in C using the code generated by the SPCIES toolbox, compiled using \texttt{gcc} version 10.2.1 with the \texttt{-O3} flag.
The \textit{equMPC} formulation uses $N=17$, which is the smallest prediction horizon for which its optimization problem was feasible during the simulation.
Note that Figures~\ref{fig:result:time} and \ref{fig:result:iter} use the right-hand-side axis to show the results for the \textit{equMPC} formulation, since they are significantly larger than the other two.
The laxMPC solver has the same number of iterations as Algorithm~\ref{alg:ADMM:MPC} but a slightly lower computation time due to its lack of terminal constraint.
However, this small computational cost comes with the advantage of better stability guarantees.

We note that we have modified the upper and lower constraints on the heat transfers to $2000$ and $-2000$, respectively, to obtain active input constraints during the first two sample times.
All other constraints can be found in \cite[Table~III]{Krupa_TCST_20}.

\section{Conclusions} \label{sec:conclusions}

We present a sparse solver for linear MPC subject to terminal quadratic constraint.
Its motivation is to be used as a substitute in applications where a polyhedral terminal set is intractable, such as in the case of large systems or in many robust MPC applications.
The proposed ADMM-based solver deals with the ellipsoidal constraint by modifying the ADMM equality constraints so that the projection step related to the ellipsoidal constraint results in an optimization problem with a simple explicit solution.
The benefit is that we retain the simple matrix structures exploited by a recently proposed MPC solver from the literature.
Our numerical results illustrate the computational advantage obtained as a result when compared against solvers which pose the quadratic constraint as a SOC constraint.
The MPC solver is suitable for its implementation in embedded systems, where its sparse nature results in a small iteration complexity and memory footprint.
In this article we focus on the case of linear MPC subject to a single (terminal) ellipsoidal constraint due to its practical usefulness.
However, the proposed method could be extended to optimization problems with multiple ellipsoidal constraints.
Finally, we remind the reader that the solver is available at \url{https://github.com/GepocUS/Spcies}.

\begin{appendices}
\renewcommand{\thesectiondis}[2]{:}

\section{Proof of Proposition \ref{theo:P:projection}} \label{app:proof:theo:P:projection}

Let $\ta \doteq a - c$, and $\Pi_{(P, c, r)}(a)$ denote the argument that minimizes \eqref{eq:P:projection:optimization:problem} for the given $P$, $c$ and $r$. It is clear that

\begin{equation} \label{eq:proof:P:proj:Pi}
    \Pi_{(P, c, r)}(a) = \Pi_{(P, 0, r)}(\ta) + c,
\end{equation}
since this simply corresponds to shifting the origin to the center of the ellipsoid $\E(P, c, r)$, and then undoing it.
Therefore, it suffices to find a closed expression for the solution of $\Pi_{(P, 0, r)}(\ta)$, i.e., to find an explicit solution to
\begin{equation}
    \min\limits_{v}\; \left\{ \frac{1}{2} \| v - \ta \|^2_P, \; {\rm s.t.}\; v\T P v \leq r^2 \right\}, \label{eq:proof:P:proj:OP}
\end{equation}
whose Lagrangian $\cc{L}: \R^n \times \R \rightarrow \R$ is given by
\begin{equation} \label{eq:proof:P:proj:Lagrangian}
    \cc{L}(v, y) = \frac{1}{2} \| v - \ta \|^2_P + y (v\T P v - r^2),
\end{equation}
where $y \in \R$ is the dual variable.
The dual function of~\eqref{eq:proof:P:proj:OP} is given by $\Psi(y) = \inf\limits_{v} \cc{L}(v, y)$.
Let us define $v(y) \inR{n}$ as $v(y) = \arg\min\limits_{v} \cc{L}(v, y)$.
It is clear from the definition of $v(y)$ that the dual function can be expressed as $\Psi(y) = \cc{L}(v(y), y)$.
The dual problem of \eqref{eq:proof:P:proj:OP} is to maximize the dual function $\Psi(y)$ subject to $y \geq 0$, i.e., to find the optimal solution of
\begin{equation} \label{eq:proof:P:proj:dual:problem}
    \max\limits_{y \geq 0}  \cc{L}(v(y), y).
\end{equation}
We start by rewritting $v(y) = \arg\min\limits_{v} \cc{L}(v, y)$ by applying simple algebraic manipulations to \eqref{eq:proof:P:proj:Lagrangian} and cancelling the terms that do not depend on $v$ and $y$, leading to
\begin{equation} \label{eq:proof:P:proj:vy}
    v(y) = \arg\min\limits_{v} \frac{1}{2} v\T (1 {+} 2 y) P v - \ta\T P v,
\end{equation}
Problem~\eqref{eq:proof:P:proj:vy} is an unconstained QP problem.
Therefore, its optimal solution is the vector $v(y)$ for which the gradient of its objective function is equal to zero \cite[\S 4.2.3]{Boyd_ConvexOptimization}.
That is, we have that $(1 {+} 2 y) P v(y) = P \ta$.
Since $y \geq 0$ is a scalar, we can divide both terms of this expression by $(1 + 2 y)$ and then multiply them by $P^{-1}$ to obtain
\begin{equation} \label{eq:proof:P:proj:v:y:solution}
    v(y) = \fracg{1}{1 {+} 2 y} \ta.
\end{equation}
Substituting expression \eqref{eq:proof:P:proj:v:y:solution} into \eqref{eq:proof:P:proj:Lagrangian} leads to
\begin{align} \label{eq:proof:P:proj:L}
    \cc{L}(v(y), y) &= \left[ \frac{1}{2} \left( \frac{2 y}{1 + 2 y} \right)^2 {+} \frac{y}{(1 + 2 y)^2} \right] \ta\T P \ta - r^2 y \nonumber \\
                    &= \fracg{y}{1 {+} 2 y} \ta\T P \ta - r^2 y,
\end{align}
which, for $y > -1/2$, is a differentiable real-valued concave function\footnote{Notice that $\fracg{y}{1 + 2 y} = \fracg{1}{2}\left( 1 - \fracg{1}{1 + 2 y} \right)$, where $-1/(1 {+} 2 y)$ is differentiable, real-valued and concave for ${y > -1/2}$.}. Therefore, given that $y$ is a scalar, the optimal solution $y^*$ of \eqref{eq:proof:P:proj:dual:problem} is given by $y^* = \max \{ \hat y, 0 \}$,
where $\hat y$ is the scalar safisfying
\begin{equation*}
    \frac{d \cc{L}(v(y), y)}{d y} \Bigr|_{\hat y} = 0.
\end{equation*}
Thus, we obtain an explicit expression for $\hat y$ by differentiating \eqref{eq:proof:P:proj:L}, resulting in
\begin{equation*}
    \frac{\ta\T P \ta - r^2 (1 {+} 2 \hat y)^2}{(1 {+} 2 \hat y)^2} = 0,\; \text{i.e.,} \;
    \hat y = \frac{1}{2} \left( \frac{\sqrt{\ta\T P \ta}}{r} - 1 \right).
\end{equation*}
Since strong duality holds, we have that the optimal solution $v^*$ is given by $v^* = v(y^*)$. Therefore, we conclude~that
\begin{align*}
    v(y^*) &= \frac{r}{\sqrt{\ta\T P \ta}} \ta,\; \text{if~} \ta\T P \ta > r^2, \\
    v(y^*) &= \ta,\; \text{if~} \ta\T P \ta \leq r^2,
\end{align*}
which, noting that $v(y^*) \equiv \Pi_{(P, 0, r)}(\ta)$ and taking into account \eqref{eq:proof:P:proj:Pi}, proves the claim. \qed

\end{appendices}

\bibliographystyle{IEEEtran}
\bibliography{IEEEabrv,BibKrupa}

\end{document}